\magnification=\magstep1
\input amstex
\documentstyle{amsppt}
\TagsOnRight
\hsize=5.1in                                                  
\vsize=7.8in

\define\C{{\Bbb C}}
\define\R{{\Bbb R}}

\define\Z{{\Bbb Z}}

\define\M{{\Cal E}}

\define\Hom{\operatorname{Hom}}

\define\Tr{\operatorname{Tr}}

\define\Det{\operatorname{Det}}

\define\GL{\operatorname{GL}}

\define\A{{\Cal A}}
\redefine\B{{\Cal B}}

\define\im{\operatorname{im}}
\define\id{\operatorname{id}}
\define\cl{\operatorname{cl}}

\define\E{\Cal E}
\define\U{{\Cal U}}
\def\<{\langle}
\def\>{\rangle}

\define\pd#1#2{\dfrac{\partial#1}{\partial#2}}
\documentstyle{amsppt}

\topmatter

\title  DETERMINANT LINES, VON NEUMANN ALGEBRAS\\
 AND L$^2$ TORSION
\endtitle
\author  A. Carey, M. Farber,  and V. Mathai\endauthor
\address
School of Mathematical Sciences,
Tel-Aviv University,
Ramat-Aviv 69978, Israel \newline (M.Farber) 
\endaddress
\address
Department of Pure Mathematics, University of Adelaide
South Australia, Australia (A.Carey and V.Mathai)
\endaddress
\vskip 2cm
\email 
acarey\@maths.adelaide.edu.au, \newline
farber\@math.tau.ac.il, \newline
vmathai\@maths.adelaide.edu.au
\endemail
\thanks{The research was supported by a grant from the
US - Israel Binational Science Foundation}
\endthanks
\abstract{In this paper, we suggest a construction
of determinant lines of finitely generated
Hilbertian modules over finite von Neumann algebras. Nonzero elements of the 
determinant lines can be viewed as volume forms on the Hilbertian modules.
Using this, we study both $L^2$ combinatorial and $L^2$ analytic 
torsion invariants
associated to flat Hilbertian bundles over compact polyhedra and manifolds;
we view them as volume forms on the reduced $L^2$ homology  and 
$L^2$ cohomology. These torsion invariants specialize to the 
the classical Reidemeister-Franz
torsion  and the Ray-Singer torsion in the finite dimensional case. Under the 
assumption that the
$L^2$ homology vanishes, the determinant line can be canonically 
identified with $\R$, and our $L^2$ torsion invariants specialize to the 
$L^2$ torsion invariants previously constructed by A.Carey, V.Mathai and
J.Lott. We also show that a recent theorem of Burghelea et al. can be reformulated
as stating equality between 
two volume forms (the combinatorial and the analytic) on the reduced 
$L^2$ cohomology.}
\endabstract
\endtopmatter

\nopagenumbers

\heading{\bf \S 0. Introduction}\endheading

The study of the $L^2$ Reidemeister-Franz torsion and
the $L^2$  analytic torsion 
was initiated in \cite{M}, \cite{L} and \cite{CM}. 
These invariants were originally defined 
for manifolds with trivial $L^2$ cohomology and positive
Novikov-Shubin invariants; they were 
shown to be piecewise linear and smooth invariants
respectively. See also \cite{LR} for K-theoretic
generalizations of these invariants.

In this paper, we introduce a new concept of 
determinant line of a finitely generated Hilbertian  
module over a von Neumann algebra. Here a Hilbertian 
module is defined as a topological vector space and a module over a von Neumann 
algebra such that there is an admissible scalar
product on it, making it a Hilbert module.
The construction of the determinant line, which we suggest here,
generalizes the classical construction of determinant line of a 
finite dimensional vector space, and 
enjoys similar functorial properties. 
Nonzero elements of the determinant line 
can be naturally
viewed as a volume forms on the Hilbertian module.
This enables us to make sense of the notions of volume forms and determinant 
lines
in the infinite dimensional and non-commutative situation.

We then define analytic and
the Reidemeister-Franz $L^2$ torsion invariants of flat Hilbertian bundles
of determinant class over finite polyhedra and compact manifolds 
respectively. These reduce to the classical constructions
in the finite dimensional situation. 
These new torsion invariants live in the determinant lines of 
reduced $L^2$ cohomology. We prove that the combinatorial $L^2$ torsion 
is combinatorially invariant (i.e. it is invariant under subdivisions).
Under the hypothesis
that the $L^2$ homology vanishes, these determinant lines can be 
canonically identified with $\R$, and our $L^2$ torsion invariants
reduce to the ones previously studied in \cite{M}, \cite{L}
and \cite{CM}.  

In a recent preprint \cite{BFKM} of Burghelea, Friedlander, Kappeler and 
McDonald, an equality  
between some combinatorial and analytic $L^2$ torsion invariants, was 
established.
More precisely, in \cite{BFKM} they introduced three numerical invariants
$T_{an}$, $T_{met}$ and $T_{comb}$ such that $T_{an}$ depends
on the Riemannian metric, $T_{comb}$ depends on the triangulation and
$T_{met}$ depends both on the Riemannian metric 
and the triangulation. The main result of \cite{BFKM} states
the equality
$T_{an} = T_{comb} \times T_{met}$.

In this paper we observe (following the suggestion of the referee) that one can 
reformulate the the main theorem of \cite{BFKM} as stating equality
between the combinatorial and analytic torsion invariants as defined
in this paper, i.e. understood as 
volume forms on the reduced $L^2$ cohomology. 

The paper is organized as follows. In the first section, which contains
preliminary material, we discuss basic properties of Hilbertian modules over
von Neumann algebras, the
canonical trace on the commutant and 
present a new construction of the Fuglede-Kadison determinant. This  
construction avoids reference to the scalar product (unlike the 
standard construction) and so it is well suited for Hilbertian modules which is 
crucial for this paper. 
In section \S 2, we construct determinant lines for finitely generated 
Hilbertian modules and establish their basic functorial properties.
In section \S 3, we consider a generalisation, which enables to deal
with $D$-admissible scalar products, $D$-isomorphisms etc.. Here 
the letter "$D$" stands for "{\it determinant class type" condition}.
In section \S 4, we define combinatorial
$L^2$ torsion, as an element in the determinant line of the reduced
$L^2$ homology and cohomology; we also prove combinatorial invariance. 
In section \S 5, we define 
$L^2$ analytic torsion as an element in the determinant line
of reduced $L^2$ cohomology.
In last section \S 6, we 
show how one may reformulate the theorem of \cite{BFKM} using the 
notions introduced in the present paper.
 
We are very thankful to the referee for his helpful comments.

\heading{\bf \S 1. Hilbertian modules and the Fuglede-Kadison determinant}
\endheading

This section contains some preliminary material which will be used 
later in this paper.

First, we describe the notion of a {\it Hilbertian module} over a finite von
Neumann algebra, which is distinct from the standard well-known 
notion of a
{\it Hilbert module}. The difference is that Hilbertian module does not 
have a specified scalar product on it, and any choice of such scalar product
introduces an additional structure. This notion will be crucial for our goals in the 
following sections. Secondly, we show that the trace on the initial
von Neumann algebra determines canonically a trace on the commutant
of any finitely generated Hilbertian module. 
At last, we use this canonical trace
to define the Fuglede-Kadison (FK) determinant. 
Note, that the commutant of a Hilbertian module is not a von Neumann
algebra (since it has no fixed $\ast$-operator). In order to be able to use the
standard theory of FK determinants, we have to show that it is independent on the
involution (which is used in its definition, cf. \cite{FK}). Instead, we provide
in this section a very simple self-contained exposition of the FK
determinants,  adjusted to needs of the current situation. 

\subheading{1.1} Let ${\A}$ be a finite  von Neumann algebra
with a fixed finite,
normal, and faithful trace $\tau:{\Cal A}\rightarrow \C$.  The involution
in ${\A}$ will be denoted $*\,$; by $\ell^{2}({\A}) = \ell^{2}_\tau({\A})$ 
we denote the
completion of ${\A}$ with respect to the scalar product $\langle a,b
\rangle=\tau(b^{*}a)$, for $a,b\in{\Cal A}$, determined by the trace $\tau$.

\subheading{1.2} Recall that a {\it Hilbert module} over ${\A}$
is a Hilbert space $M$
together with a continuous left ${\Cal A}$-module structure such that there
exists an isometric ${\Cal A}$-linear embedding of $M$ into
$\ell^{2}({\Cal A})\otimes H$, for some Hilbert space $H$.
Note that this embedding is not
part of the structure.  A Hilbert module $M$ is  said to be
{\it finitely generated} if it
admits an embedding $M\rightarrow\ell^{2}({\Cal A})\otimes H$ as above with
finite dimensional $H$.

Any Hilbert module, being a Hilbert space, has a particular scalar product.
In this paper we
wish to consider a weaker notion obtained from Hilbert module by forgetting
the scalar product but preserving its topology and the ${\Cal A}$-action.

\subheading{1.3. Definition} A {\it  Hilbertian module}
is a topological vector space $M$
with continuous left ${\Cal A}$-action such that there exists a scalar product
$\langle\;,\;\rangle$ on $M$ which generates the topology of $M$ and such
that $M$ together
with $\langle\;,\;\rangle$ and with the ${\Cal A}$-action is a Hilbert
module.

If $M$ is a Hilbertian module,
then any scalar product $\langle\;,\;\rangle$ on
$M$ with the above properties will be called {\it admissible}.

A {\it morphism} of Hilbertian modules  $f:M\to N$ is a continuous linear map
commuting with the $\A$-action.

Note that the kernel of any morphism $f$ is again a Hilbertian module. Also,
the closure of the image $\cl(\im(f))$ is a Hilbertian module.

\subheading{1.4} Let $\langle\;,\;\rangle$ be an admissible scalar product on $M$. 
Then $\langle\;,\;\rangle$ must be compatible with the
topology on $M$ and with the ${\Cal A}$-action.
The last condition means that
the involution on ${\Cal A}$ determined by the scalar product $\langle\;,
\;\rangle$ coincides with the involution of the von Neumann algebra
${\Cal A}:$
$$
\langle\lambda\cdot v,w\rangle=\langle v,\lambda^{*}\cdot w\rangle\tag1
$$
for any $v,w\in M$, $\lambda\in{\Cal A}$.

Suppose now that $\langle\;,\;\rangle_1$ is another admissible
scalar product.  Then there exists an operator
$$
A:M\rightarrow M
$$
such that
$$
\langle v,w\rangle_1=\langle Av,w\rangle\tag2
$$
for any $v,w\in M$.  The operator $A$ has to be:
\roster
\item"($\alpha$)"  {\it a linear homeomorphism 
(since the scalar products $\<\ ,\ \>$ and
$\<\ ,\ \>_1$ define the same topology);
\item"($\beta$)" self-adjoint;
\item"($\gamma$)" positive;
\item"($\delta$)" commuting with the $\A$-action.}
\endroster

In order to prove the last property, write
$$
\align
\<\lambda\cdot v,w\>_1 &= \<v,\lambda^*\cdot w\>_1,\\
\<A(\lambda\cdot v),w\> &= \<A(v),\lambda^*\cdot w\>,\\
\<A(\lambda\cdot v),w\> &= \<\lambda\cdot A(v), w\>
\endalign
$$
for any $v,w\in M,\quad \lambda\in \A$. Thus, it follows that
$A(\lambda v)\ =\ \lambda A(v).$

We conclude: {\it once an admissible scalar product on a Hilbertian
module $M$ has been chosen, there is a one-to-one correspondence between the
admissible scalar products on $M$ and the operators $A:M\to M$ satisfying
($\alpha$) - ($\delta$) above.}

\proclaim{1.5. Corollary} If $\<\ ,\ \>$ and $\<\ ,\ \>_1$ are two admissible
scalar products on a Hilbertian module $M$  then the Hilbert modules
$(M, \<\ ,\ \>)$ and $(M, \<\ ,\ \>_1)$ are isomorphic.
\endproclaim
\demo{Proof} Let $A:M\to M$ be the operator such that
$\<v,w\>_1\ =\ \<Av,w\>$ for $v,w\in M$.
Let $B:M\to M$ be the positive square root of $A$; since the spectrum of
$A$ is real and positive, the functional
calculus produces the operator $B$ uniquely. The operator $B$ is positive,
invertible, and commutes with the action of $\A$. Hence
the map $x\mapsto Bx,\ x\in M$
establishes an isomorphism $(M, \<\ ,\ \>_1)\to (M, \<\ ,\ \>)$.
$\square$
\enddemo

\subheading{1.6} Using this corollary we may define
{\it finitely generated Hilbertian
modules} as those for which the corresponding Hilbert modules
(obtained by a choice of an admissible scalar product) are finitely
generated.

Note that the {\it von Neumann dimension}
(denoted $\dim_{\tau}(M)$ or $\tau(M)$)
of a Hilbertian module is also correctly defined (by virtue of
the Corollary 1.5 above).

\subheading{1.7. The commutant} Let $M$ be a Hilbertian ${\Cal A}$-module.
Let $\B\ =\ \B(M)$ denote the algebra $\B=B_{\Cal A}(M)$
of all bounded
linear operators on $M$ commuting with ${\Cal A}$ (the {\it commutant of}
$M$).

Any choice of an admissible scalar product $\langle\;,\;\rangle$ on $M$,
defines obviously a $*$-operator on $\B$ (by assigning to an operator
its adjoint); this turns $\B$ into a von Neumann algebra.  Note that this
involution $*$ depends on the scalar product $\<\ ,\ \>$ on $M$;
if we choose another admissible scalar product $\<\ ,\ \>_1$
on $M$, then the new involution will be given by
$$
f\mapsto \ A^{-1}f^*A\qquad \text{for}\quad f\in \B,\tag3
$$
where $A\in\B$ is the operator defined by
$\<v,w\>_1\ =\ \<Av,w\>$ for $v,w\in M$.

Our aim now is to construct canonically a trace on $\B$, using the given
trace $\tau$ on ${\Cal A}$.

\proclaim{1.8. Proposition} If $M$ is a finitely generated Hilbertian module, then
the trace $\tau:{\Cal A}\rightarrow \C$
determines canonically a trace on the commutant
$$
\Tr_{\tau}:\B\ =\ \B(M)\rightarrow \C\tag4
$$
which is finite, normal, and faithful (with respect to the involution on $\B$ 
corresponding to any choice of an admissible scalar product on $M$).
If $M$ and $N$ are two finitely
generated modules over $\A$, then the canonical traces $\Tr_\tau$ on
$\B(M)$, $\B(N)$ and on $\B(M\oplus N)$ are compatible 
in the following sense:
$$\Tr_\tau \left(\matrix
A&B\\
C&D
\endmatrix \right)
= \Tr_\tau (A) \ +\ \Tr_\tau(D),\tag5$$
for all $A\in \B(M),\quad D\in\B(N)$ and any 
morphisms $B:M\to N$, and $C:N\to M$.
\endproclaim

\demo{Proof} The proof is straightforward and therefore we will indicate only
the main steps. 

Suppose first that $M$ is {\it free},
that is, $M$ is isomorphic
to $l^2(\A)\otimes\C^k$ for some $k$. Then the commutant
$\B(M)$ can be identified with the algebra of
$k\times k$-matrices with entries in $\A$, acting from the right on
$l^2(\A)\otimes\C^k$ (the last module is viewed as the set of row-vectors
with entries in $l^2(\A)$). If $\alpha\in\B$ is represented
by a $k\times k$ matrix $(\alpha_{ij}),$ then one defines
$$ \Tr_{\tau}(\alpha)\ =\ \sum_{i=1}^k \ \tau(\alpha_{ii}).$$
This gives a trace on $\B$ which satisfies all necessary conditions.

If $M$ is not free, then we can embed it in a free module as a closed 
$\A$-invariant subspace. Then the commutant $\B(M)$ can be identified with 
a left ideal in the $k\times k$-matrix algebra with entries in $\A$ and 
the trace described in the previous paragraph restricts to this ideal and 
determines a trace on $\B(M)$. One then shows that the obtained 
trace on $\B(M)$ does not depend on the embedding of $M$ in a free module.
\qed
\enddemo

\subheading{1.9} Let $M$ be a finitely generated Hilbertian module
over a von Neumann algebra $\A$, as above. Denote by $\GL(M)$ the group of
all invertible elements
of the commutant $\B(M)=B_{\A}(M)$ . We will consider
the norm topology on $\GL(M)$; with this topology it is a Banach Lie group.
Its Lie algebra can be identified with the commutant $\B(M)$. The canonical
trace $\Tr_{\tau}$ on the commutant $\B(M)$ (described in the previous
subsection) is a homomorphism of the Lie algebra $\B(M)$ into the abelian
Lie algebra $\C$. By the standard theorems, it defines a group homomorphism
of the universal covering group of $\GL(M)$ into $\C$. This approach
leads to following construction of the Fuglede-Kadison determinants,
compare \cite{HS}.

\proclaim{1.10. Theorem} There exists a function
$$\Det_{\tau}: \GL(M)\to \R^{>0}\tag7$$
(called the Fuglede - Kadison determinant), satisfying:
\roster
\item"(a)" $\Det_{\tau}$ is a group homomorphism, that is,
$$\Det_{\tau}(AB)\ =\ \Det_{\tau}(A)\cdot \Det_{\tau}(B)\tag8$$
for $A,\ B\in\GL(M)$;
\item "(b)" $$\Det_{\tau}(\lambda I)\ =\ |\lambda|^{\dim_{\A}(M)}\tag9$$
for $\lambda \in \C$, $\lambda\ne 0$; here $I\in\GL(M)$ denotes the identity
operator;
\item "(c)" $$\Det_{\lambda\tau}(A) \ =\ \Det_{\tau}(A)^\lambda\tag10$$
for $\lambda\in \R^{>0}$;
\item"(d)"  $\Det_{\tau}$ is continuous as a map $\GL(M)\to \R^{>0}$, where
$\GL(M)$ is supplied with the norm topology;
\item "(e)" If $A_t$ for $t\in [0,1]$
is a continuous piecewise smooth path in $\GL(M)$ then
$$
\log \lbrack\frac{\Det_{\tau}(A_1)}{\Det_{\tau}(A_0)}\rbrack\ =\
\int_0^1\Re\Tr_{\tau}\lbrack A_t^{-1}A_t^\prime\rbrack dt.\tag11
$$
Here $\Re$ denotes the real part, $\Tr_{\tau}$ denotes the canonical
trace on the commutant constructed in Proposition 1.8 and $A_t^\prime$ 
denotes the derivative of $A_t$ with respect to $t$.
\item"(f)"  Let $M$ and $N$ be two finitely generated modules over $\A$, and
$A\in\GL(M)$ and $B\in\GL(N)$ two invertible automorphisms, and let 
$\gamma:N\to M$ be a homomorphism. Then the map given by the matrix
$$\left(\matrix
A&\gamma\\
0&B
\endmatrix\right)$$
belongs to $\GL(M\oplus N)$ and
$$
\Det_\tau\left(\matrix
A&\gamma\\
0&B
\endmatrix\right)\ =\
\Det_\tau(A)\cdot\Det_\tau(B)\tag12
$$
\endroster

\endproclaim

\demo{Proof} We are going to accept the following definition.
Given an invertible operator $A\in\GL(M)$, find a continuous piecewise
smooth path
$A_t\in\GL(M)$ with $t\in [0,1]$, such that $A_0=I$ and $A_1=A$ (it is well
known that the group $\GL(M)$ is pathwise connected, cf. \cite{Di}). Then define
$$\log \Det_{\tau}(A)\ =\
\int_0^1\Re\Tr_{\tau}\lbrack A_t^{-1}A_t^\prime\rbrack dt.\tag13$$

We want to show that the integral does not depend on the choice of the path,
joining $A$ with the identity $I$.
Consider the integrals of the form
$$\int_a^b \Tr_{\tau}[A^{-1}_tA^\prime_t]dt,\tag14$$
where $A_t$ for $t\in [a,b]$ is a piecewise smooth path in $\GL(M)$.
Suppose first that the path $A_t$ is "small" in the
following sense: for all $t\in [a,b]$ we have $||A_a^{-1}A_t-I||<1$.
Then one may find a piecewise smooth path $C_t\in\B(M)$ such that
$A_t=A_a\exp(C_t)$ and $C_a=0$; here $C_t$ is defined as
$\log (A_a^{-1}A_t)$ and the
logarithm function defined by its Taylor power expansion around 1.
Using Duhamel's formula
$$A^\prime_t\ =\ A_a\int_0^1 e^{(1-s)C_t}C_t^\prime e^{sC_t} ds\tag15$$
one computes
$$
\align
\int_a^b \Tr_{\tau}[A_t^{-1}A_t^\prime]dt
&= \int_a^b\Tr_{\tau}(C_t^\prime)dt \ =\\
&= \int_a^b\frac{d}{dt}(\Tr_{\tau}(C_t))dt\ =\ \Tr_{\tau}(C_b) \ =\\
&= \Tr_{\tau}(\log(A_a^{-1}A_b))\tag16
\endalign
$$

If $A_t$ is an arbitrary piecewise smooth path in $\GL(M)$, defined for
$t\in [0,1]$, (which is now not supposed to be "small"), then the interval
$[0,1]$ can be divided into subintervals $t_0=0<t_1<\dots <t_N=1$ such that
$||A_{t_i}^{-1}A_t-I||<1$ for all $t\in [t_i,t_{i+1}]$. Thus by the argument
above we have
$$
\int_a^b \Tr_{\tau}[A_t^{-1}A_t^\prime]dt \ =\
\sum_{i=0}^{N-1}\Tr_{\tau}(\log(A_{t_i}^{-1}A_{t_{i+1}}))\tag17
$$
This shows, in particular, that the integral (14)
depends only on the homotopy class of the path (relative to the end points).
One easily checks the following homomorphism property: if $A_t$ and $B_t$
are two piecewise smooth paths in $\GL(M)$ defined for $t\in [a,b]$ and if $C_t=A_tB_t$
is their product then
$$\int_a^b \Tr_{\tau}[C_t^{-1}C_t^\prime]dt \ =\
\int_a^b \Tr_{\tau}[A_t^{-1}A_t^\prime]dt\ +\
\int_a^b \Tr_{\tau}[B_t^{-1}B_t^\prime]dt.\tag18$$

Now we observe that the integral
$$\int_a^b \Re\Tr_{\tau}[A_t^{-1}A_t^\prime]dt,\tag19 $$
taken along any closed loop in $\GL(M)$, vanishes. This follows from the
fact that any closed loop in $\GL(M)$ is homotopic
to a product of loops of the form
$$A_t\ =\ \exp(2\pi itS),\qquad\text{for}\quad 0\le t\le 1,\tag20$$
where $\in\B(M)$ and $\Tr_\tau(S)$ is real; for loops of this form the
integral vanishes obviously. In fact, if we choose an
admissible scalar product on $M$ then the commutant $\B(M)$ becomes a
finite von Neumann algebra; the result of H. Araki, M-S.B. Smith and L. Smith
(cf. \cite{ASS}, Theorem 2.8) states that any loop in the group $\U(M)$ of
unitary elements of $\B(M)$ is homotopic to a product of the loops of the
form (20)
where $S$ is a self-adjoint element of $\B(M)$. Using the polar decomposition
it is easy to construct a deformation retraction of $\GL(M)$ onto
$\U(M)$ and thus the result follows. 

This gives a construction of the FK determinant $\Det_\tau$. The homomorphism
property (18) proves multiplicativity of the determinants (a). The other
properties (b), (c), (d), (e) are clearly satisfied. We are left to prove (f).
Let $A_t$ and $B_t$ be piecewise smooth paths connecting $A$ and $B$ with
the identity maps of $M$ and $N$ correspondingly. Then the family of maps
$\left(\matrix A_t&t\gamma\\0&B_t\endmatrix \right)$ joins the given map 
with the identity of $M\oplus N$ and we obtain
$$
\aligned
&\Tr_\tau\left( \left(
\matrix A_t&t\gamma\\ 0&B_t\endmatrix\right)^{-1}\cdot
\left(\matrix A_t^\prime &\gamma\\ 0&B_t^\prime\endmatrix\right) \right)\ = \\
&\Tr_\tau\left( \left(
\matrix A_t^{-1}& -B_t^{-1}t\gamma A_t^{-1}\\
0& B_t^{-1}\endmatrix\right)\cdot\left(\matrix A^\prime_t&\gamma\\
0&B^\prime_t\endmatrix\right) \right)\ =\\ 
&\Tr_\tau\left(\matrix A^\prime_tA_t^{-1}& \ast \\
0&B^\prime_tB_t^{-1}\endmatrix\right) \ =\\ 
&\Tr_\tau(A^\prime_tA_t^{-1})\ +\ 
\Tr_\tau(B^\prime_tB_t^{-1})
\endaligned
$$
where the formula (5) of Proposition 1.8 has been used. 
\qed
\enddemo

\subheading{1.11} As an example of computation of the Fuglede-Kadison 
determinant
based on the above definition, consider the following situation.
Suppose that an operator $A\in \GL(M)$ is given by
$$A\ =\ \int_0^\infty \lambda dE_\lambda\tag21$$
where $\lambda$ is real and $E_\lambda$ is a spectral projection in $\B(M)$. 
We will assume that the operator $A$ is
invertible in $\B(M)$; this means that $E_\lambda=0$ for sufficiently
small $\lambda$. Then we can choose the straight path
$$A_t\ =\ t(A-I)\ +\ I,\quad t\in [0,1]$$
as the path joining $A$ with $I$ inside $\GL(M)$. Applying the above
definition we obtain
$$
\align
&\log \Det_\tau(A)\ = \int_0^1\Re\Tr_\tau[A_t^{-1}A_t^\prime]dt = \\
&\int_0^\infty [\int_0^1\frac{\lambda -1}{t(\lambda-1)+1}dt]d\phi_\lambda
= \int_0^\infty \log\lambda d\phi_\lambda,\tag22
\endalign
$$
where $\phi_\lambda = \Tr_\tau E_\lambda$ is the spectral density function.

\heading{\bf \S2.  Determinant line of a Hilbertian module}\endheading

In this section we describe a construction of the determinant line associated to
a finitely generated Hilbertian module. It will be used later in \S 3
in the constructions of combinatorial and analytic torsion invariants 
associated to polyhedra and compact manifolds.

Let $\A$ be a finite von Neumann algebra with a fixed finite, normal, and
faithful trace $\tau: \A\to \C$.
Given a Hilbertian module $M$ of finite type over $\A$,
we are going to associate with a $M$ (in a canonical way) an oriented real
line which we will denote $\det(M)$; we will call it the {\it determinant
line of} $M$. Our construction will generalize the determinant line of
a finite dimensional vector space.

In order to clarify the precise definition (given in the paragraph below)
let's make the following remark. In the category of finite dimensional
vector spaces, the determinant line consists of volume forms. What is the
correct generalization of the notion of volume form for Hilbertian modules?
Note that one can define a volume form on a vector space by presenting
a scalar product on this vector space and two different scalar products
determine the same volume form if and only if the determinant of the transition
operator (determined by the pair of scalar products) is equal to 1.
Thus, we may consider a volume form as an equivalence class of admissible
scalar products. In fact, in the last form, the notion of volume form
can be generalized to Hilbertian modules over von Neumann algebras,
using the FK determinants.

\subheading{2.1} Define $\det(M)$ 
as a real vector
space generated by symbols
$\<\ ,\ \>$, one for any admissible scalar product on $M$, subject to the
following relations: for any pair  $\<\ ,\ \>_1$ and $\<\ ,\ \>_2$
of admissible scalar products on $M$ we write the following relation
$$\<\ ,\ \>_2\ =\ {\Det_\tau(A)}^{-1/2} \cdot \<\ ,\ \>_1,\tag23$$
where $A\in \GL(M)$ is such that
$$\<v,w\>_2\ =\ \<Av,w\>_1$$
for all $v,w\in M$. Here the transition operator $A$ is invertible and
belongs to the
commutant $\B(M)$ (cf. section 1) and $\Det_\tau(A)$ denotes the
Fuglede-Kadison determinant of $A$ constructed (in subsection 1.10)
with the aid of the canonical trace $\Tr_\tau$ on the commutant.

Assume that we have three different admissible scalar products
$\<\ ,\ \>_1$, $\<\ ,\ \>_2$ and $\<\ ,\ \>_3$ on $M$. Suppose that
$$\<v,w\>_2=\<Av,w\>_1\quad\text{and}\quad \<v,w\>_3=\<Bv,w\>_2$$
for all $v,w\in\B(M)$, where $A,B\in \B(M)$, then
$\<v,w\>_3=\<ABv,w\>_1$ and in $\det(M)$ we have the relations
$$\<\ ,\ \>_2\ =\ {\Det_\tau(A)}^{-1/2} \cdot \<\ ,\ \>_1,$$
$$\<\ ,\ \>_3\ =\ {\Det_\tau(B)}^{-1/2} \cdot \<\ ,\ \>_2,$$
$$\<\ ,\ \>_3\ =\ {\Det_\tau(AB)}^{-1/2} \cdot \<\ ,\ \>_1,$$
and the third relation follows from the first two via the homomorphism
property of the Fuglede-Kadison determinant.

Thus we obtain, that $\det(M)$ {\it is one-dimensional real vector space
generated by the
symbol $\<\ ,\ \>$ of any admissible scalar product on $M$}.

Note also, that the real line $\det(M)$ has {\it the canonical orientation}, 
since the
transition coefficients ${\Det_\tau(A)}^{-1/2}$ are always positive. Thus we
may speak of {\it positive and negative} elements of $\det(M)$. The set of all 
positive elements of $\det(M)$ will be denoted $\det_+(M)$.

We will think of elements of $\det(M)$ as "volume forms" on $M$.

If $M$ is trivial module, $M=0$, then we set $\det(M)=\R$, by definition.

To illustrate our definitions, consider the case $\A = \C$ with the standard 
trace.
Now $M$ is just a finite dimensional vector space over $\C$. Any scalar 
product $<\ ,\ >$ on $M$ determines an element
$$e_1\wedge e_2\wedge  \dots \wedge e_n \in \Lambda^{n}_{\C}(M).$$
Here $n = \dim_{\C}M$, $\Lambda^{n}_{\C}(M)$ denotes the highest exterior 
power of $M$ and 
$e_1, e_2, \dots, e_n$ is an orthonormal basis of $M$ with respect to
$<\ ,\ >$. We obtain a map 
$\det_+(M) \to |\Lambda_{\C}^n(M)| = \Lambda^{n}(M)/\sim$; the equivalence 
relation $\sim$ is $v\sim w$ iff $v=e^{i\phi} w$ for some $\phi\in \R$.
This map is well defined; to show this one observes that 
in finite dimenesinal case the Fuglede - Kadison
determinant coincides with the absolute value of the usual determinant.
Also, if we multiply the 
scalar product $<\ ,\ >$ by $\lambda^2$, where $\lambda>0$,
then the new orthonormal basis will be
$\lambda^{-1}e_1,$ $\lambda^{-1}e_2,$ $\dots, \lambda^{-1}e_n$ and the 
corresponding element of $\Lambda^{n}_{\C}(M)$ will be 
$\lambda^{-n}e_1\wedge e_2\wedge \dots \wedge e_n$; this is compatible
with (23) and shows that the map  
$\det_+(M) \to |\Lambda^n_{\C}(M)|$ is $\R_+$-linear.
Thus, we obtain that in finite dimensions we have an identification
$\det_+(M)\simeq |\Lambda^n_{\C}(M)|$.

\subheading{2.2} Given two finitely generated Hilbertian modules
$M$ and $N$ over $\A$,
and a pair $\<\ ,\ \>_M$ and $\<\ ,\ \>_N$ of admissible scalar
products on $M$ and $N$ correspondingly, we may obviously define the
scalar product $\<\ ,\ \>_M \oplus \<\ ,\ \>_N$ on the direct sum
 $M\oplus N$. This
defines an isomorphism
$$\det(M)\otimes\det(N)\to\det(M\oplus N).\tag24$$
Using Theorem 1.10, it is easy to show that
{\it this homomorphism is canonical}, that is, it does not depend on the
choice of the metrics
$\<\ ,\ \>_M$ and $\<\ ,\ \>_N$. From the description given above it is
clear that the homomorphism (22) preserves the orientations.

\subheading{2.3} Note that, {\it any isomorphism $f:M\to N$ between
finitely generated
Hilbertian modules induces canonically an isomorphism of the determinant
lines
$$f_\ast:\det(M)\to\det(N).\tag25$$
Moreover, the induced map $f_\ast$ preserves the orientations of the
determinant lines.}
Indeed, if $\<\ ,\ \>_M$ is an admissible scalar product on
$M$, then we set
$$f_\ast(\<\ ,\ \>_M)= \<\ ,\ \>_N,\tag26$$
where $\<\ ,\ \>_N$ is the scalar product on $N$ given by
$\<v,w\>_N=\<f^{-1}(v),f^{-1}(w)\>_M$ for $v,w\in N$ 
(this scalar product it admissible
since $f$ is an isomorphism). This definition does not depend on the choice
of the scalar product $\<\ ,\ \>_M$ on $M$: if we have a different admissible
scalar product $\<\ ,\ \>_M^\prime$ on $M$, where $\<v,w\>_M^\prime =
\<A(v),w\>_M$ with $A\in\GL(M)$, then the induced scalar product on $N$
will be
$$\<v,w\>_N^\prime\ =\ \<({f}^{-1}Af)v,w\>_N$$
and our statement follows from property (a) of the Fuglede-Kadison
determinant, cf. Theorem 1.10.

\subheading{2.4} It is obvious from the definition, that our construction
is {\it functorial}:
if $f:M\to N$ and $g:N\to L$ are two isomorphisms between finitely generated
Hilbertian modules then
$$(g\circ f)_\ast\ =\  g_\ast\circ f_\ast.\tag27$$

\proclaim{2.5. Proposition} If $f:M\to M$ is an automorphism of a
finitely generated Hilbertian module $M$, $f\in\GL(M)$, then the induced
homomorphism $f_\ast:\det(M)\to \det(M)$ coincides with the multiplication
by $\Det_\tau(f)\in \R^{>0}$.
\endproclaim
\demo{Proof} Let $\<\ ,\ \>$ be an admissible scalar product on $M$, then
the induced by $f$ scalar product $\<\ ,\ \>^\prime$ is given by
$$\<v,w\>^\prime\ =\ \<f^{-1}(v),f^{-1}(w)\> \ = \ \<(ff^\ast)^{-1}v,w\>.$$
Thus in $\det(M)$ we have
$$\<\ ,\ \>^\prime \ =\ \sqrt{\Det_\tau(ff^\ast)}\cdot\<\ ,\ \>\ =
\Det_\tau(f)\cdot\<\ ,\ \>.$$
\enddemo

\proclaim{2.6. Proposition} Any exact sequence
$$0\to M^\prime@>{\alpha}>>M@>{\beta}>>M^{\prime\prime}\to 0\tag28$$
of finitely generated Hilbertian modules determines canonically an isomorphism
$$\det(M^\prime)\otimes\det(M^{\prime\prime})\ \to \det(M),\tag29$$
preserving the orientations.
\endproclaim
\demo{Proof} As is well known, any such exact sequence splits. Choose a splitting
$$0\to M^\prime@<r<< M@<s<< M^{\prime\prime}\to 0.$$
Then any pair of admissible scalar products $\<\ ,\ \>_{M^\prime}$ and
$\<\ ,\ \>_{M^{\prime\prime}}$ on $M^\prime$ and $M^{\prime\prime}$
correspondingly, define the following admissible scalar product on $M$:
$$\<v,w\>_M\ =\ \<r(v),r(w)\>_{M^\prime}\ +\
\<\beta(v),\beta(w)\>_{M^{\prime\prime}}\tag30$$
This defines an isomorphism
$$\det(M^\prime)\otimes\det(M^{\prime\prime})\ \to \det(M).$$
This isomorphism can be also described as the composition
$$\det(M^\prime)\otimes\det(M^{\prime\prime})\to \det(M^\prime\oplus
M^{\prime\prime})@>{(\alpha\oplus s)_\ast}>>\det(M)$$
where the first map is given by (24). To finish the proof we only have 
to show independence on
the splitting $s$. If $s^\prime$ is another splitting then it can be
represented in the form
$$s^\prime =\ s\ +\ \alpha\circ\gamma,\tag31$$
where $\gamma:M^{\prime\prime}\to M^\prime$ is a homomorphism.
Now the result follows from the commutative diagram
$$
\CD
M^\prime\oplus M^{\prime\prime} @>{\alpha\oplus s}>> M\\
@A
{\left(\smallmatrix 1&\gamma\\
0&1\endsmallmatrix\right)}AA  @AA{=}A   \\
M^\prime\oplus M^{\prime\prime} @>>{\alpha\oplus s'}> M
\endCD\tag32
$$
and from the fact that $\Det_\tau
\left(\smallmatrix 1&\gamma\\
0&1\endsmallmatrix\right) = \ 1,
$ cf. Theorem 1.10, statement (f). \qed
\enddemo

\heading{\bf \S3. $D$-admissible scalar products, $D$-isomorphisms, and
$D$-exact sequences}\endheading

In this section, we consider generalisations of the notions and results of section 2.
We introduce the notion of a $D$-admissible scalar product on a finitely 
generated Hilbertian module (here the letter $D$ refers to the {\it "determinant
class" type condition}, following the terminological convention suggested 
in \cite{BFKM}). We show that any $D$-admissible scalar product defines a 
nonzero element 
in the determinant line. We also introduce the notion of $D$-isomorphism and
$D$-exact sequences. We prove that any $D$-isomorphism induces an 
isomorphism of the
determinant lines. Under the hypothesis that the chain 
complex is of determinant class, we construct a natural isomorphism between the
determinant line of the chain complex and the determinant line of its 
reduced $L^2$ homology.

\subheading{3.1} Consider a Hilbertian module $M$ of finite type over $\A$ 
(as above).
A scalar product $\<\ ,\ \>$ on $M$ will be called {\it $D$-admissible}
if it can be represented in the form
$$\<v,w\>\ =\ \<A(v),w\>_1\quad\text{for}\quad v,w\in M,$$
where $\<\ ,\ \>_1$ is an admissible scalar product on $M$ and $A\in\B_{\A}(M)$
is an injective (possibly not invertible) homomorphism $A:M\to M$, which
is positive and self-adjoint with respect to $\<\ ,\ \>_1$, and
the following property is satisfied: if
$$A\ =\ \int_0^\infty \lambda dE_\lambda\tag33$$
is the spectral decomposition of $A$ and if
$$\phi(\lambda)\ =\ \dim_\tau(E_\lambda)\ =\ \Tr_\tau(E_\lambda)\tag34$$
denotes the corresponding spectral density function, then the integral
$$\int_0^\infty\ln(\lambda)d\phi(\lambda)\ >\ -\infty\tag35$$
is assumed to converge to a finite value. Note that the integral (35) may only
diverge at point $\lambda=0$. 

We want to show that {\it any $D$-admissible scalar product determines
canonically a nonzero element of the determinant line $\det(M)$}.

\proclaim{3.2. Proposition} In the above notations, the non-zero element
$$\exp[-1/2\int_0^\infty\ln(\lambda)d\phi(\lambda)]\cdot \<\ ,\ \>_1\tag36$$
of the determinant line $\det(M)$ depends only on the $D$-admissible
scalar product $\<\ ,\ \>$ and does not depend on the choice of the scalar
product $\<\ ,\ \>_1$. The class of this element in $\det(M)$ will be
denoted by the same symbol $\<\ ,\ \>$.
\endproclaim

Introduce the following notation
$$\Det_\tau(A)=\exp[\int_0^\infty\ln(\lambda)d\phi(\lambda)];\tag37$$
it is an extension of the Fuglede-Kadison determinant to some non-invertible
self-adjoint positive operators, which was discussed
in \cite{FK},\cite{Lu} (compare (22)).
With this notation the formula above can be written as
$$\<\ ,\ \>\ =\ \Det_\tau(A)^{-1/2}\cdot \<\ ,\ \>_1,\tag38$$
which is consistent with our previous constructions.

\demo{Proof of Proposition 3.2} Suppose that $\<\ ,\ \>_2$ is another
admissible scalar
product on $M$ and let $\<v,w\>=\<Bv,w\>_2$ for $v,w\in M$,
where $B\in\B(M)$. Then there is an invertible $C\in\GL(M)$,
such that $\<v,w\>_2=\<Cv,w\>_1$, where $C\in\GL(M)$. Hence $A=CB$ and
we obtain
$$\Det_\tau(B)^{-1/2}\cdot\<\ ,\ \>_2\ =\
\Det_\tau(B)^{-1/2}\cdot\Det_\tau(C)^{-1/2}\cdot \<\ ,\ \>_1=\
\Det_\tau(A)^{-1/2}\cdot\<\ ,\ \>_1$$
since $\Det_\tau(A)=\Det_\tau(C)\Det_\tau(B)$, cf. \cite{FK}, \cite{Lu}.
$\square$
\enddemo

\subheading{3.3. Definition} An injective with dense image homomorphism 
$f:M\to N$ of Hilbertian modules will be called
{\it $D$-isomorphism} if for some (and hence for any) admissible
scalar product $\<\ ,\ \>_N$ on $N$ the induced scalar product
$\<\ ,\ \>_M$ on $M$ (which is given by
$\<v,w\>_M=\<f(v),f(w)\>_N$ ) is $D$-admissible.

Any $D$-isomorphism $f:M\to N$ induces an isomorphism
$$f_\ast:\det(M)\to\det(N),$$
where $f^\ast(\<\ ,\ \>_M)\ =\ \<\ ,\ \>_N$ in the above notations. 

\subheading{3.4. Definition} A sequence of Hilbertian modules and homomorphisms
$$0\to M^\prime@>{\alpha}>>M@>{\beta}>>M^{\prime\prime}\to 0\tag39$$
will be called {\it $D$-exact} if $\alpha$ is a monomorphism,
$\im(\alpha)\ =\ \ker(\beta)$, and the induced by $\beta$ map
$M/\ker(\beta)\to M^{\prime\prime}$ is a $D$-isomorphism in the sense
of Definition 3.3.

We can slightly generalize Proposition 2.6:

\proclaim{3.5. Proposition} Any $D$-exact sequence
$$0\to M^\prime@>{\alpha}>>M@>{\beta}>>M^{\prime\prime}\to 0$$
determines canonically an isomorphism
$$\det(M^\prime)\otimes\det(M^{\prime\prime})\ \to \det(M).\tag40$$
\endproclaim
\demo{Proof} Denote $\bar M\ =\ M/\ker(\beta)$ and let $\bar\beta:\bar M\to
M^{\prime\prime}$ be the homomorphism determined by $\beta$. It is a 
$D$-isomorphism and so by the remark above it induced an
isomorphism
$${\bar\beta}_\ast: \det(\bar M)\to \det(M^{\prime\prime}).$$
Thus we obtain the following isomorphism
$$\det(M^\prime)\otimes\det(M^{\prime\prime})
@>{\id\otimes\bar\beta_\ast^{-1}}>>
\det(M^\prime)\otimes\det(\bar M)@>>>\det(M),$$
where the last isomorphism is given by proposition 2.6, applied to
the exact sequence
$$0\to M^\prime@>{\alpha}>>M@>>>\bar M\to 0.$$
This completes the proof.
\enddemo

\subheading{3.6} Let $M_\ast\ =\ \oplus M_i$ be a
{\it graded Hilbertian $\A$-module of
finite type}. This means that each $M_i$ is finitely generated and there
are only finitely many nonzero modules $M_i$. We define the determinant line
of $M$ in the usual way:
$$\det(M)\ =\ \otimes \det(M_i)^{(-1)^i},\tag41$$
where $\det(M_i)^{-1}$ denotes the dual line of $\det(M_i)$.

\subheading{3.7} A crucial role
in the finite dimensional linear algebra of determinant lines is played
by the canonical isomorphism between the determinant line of a chain complex
and the determinant line of its homology, cf. \cite{BGS}. The similar
statement in the category of Hilbertian modules {\it is not true in general},
but it is true under an additional requirement, which we are
going to describe now. 
Roughly speaking, this condition means that {\it the reduced $L^2$ cohomology
"properly represents" the "full" cohomology}. Recall that the reduced 
$L^2$ cohomology is the factor of the submodule of cycles by the {\it closure} 
of the submodule of the boundaries, cf. (45).

\subheading{3.8. Definition} Let
$$0\to C_N@>{\partial}>> C_{N-1}@>{\partial}>> \dots C_0\to 0\tag42$$
be a chain complex of finite length formed by finitely generated Hilbertian
modules and bounded linear maps $\partial$ commuting with the $\A$-action.
Let $Z_i$ denote the submodule of cycles and let $B_i$ denote the submodule
of boundaries. Following \cite{BFKM}, we will say that the chain complex (42) 
belongs to  {\it the determinant class} if the following sequence
$$0\to Z_i\to C_i\to \bar B_{i-1}\to 0\tag43$$
is {\it $D$-exact}, cf. 3.4 above. 

Note, that the above condition {\it imposes no restrictions on
the size of the reduced $L^2$-homology} $Z_i/\bar B_i$, 
which may be arbitrary large. This condition can be expressed through the 
density functions
of the Laplacians of different dimensions; it is satisfied if the Novikov-Shubin
invariants are all positive. 

More precisely, the property of chain complex (42) to
be of determinant class {\it depends only on the torsion part of its  extended 
$L^2$ homology}, as defined in \cite{F1}, \cite{F2}. In fact, the torsion part of the
extended $L^2$ homology of complex (42)  in dimension $i-1$  equals to 
$(\partial: C_i\to\overline B_{i-1})$, cf.  formulae (16) and (18) in 
\cite{F2}, and our statement follows by comparing the Definitions 3.1, 3.3,
3.4, and 3.8.
The proof of the fact that the property (35) depends only on the {\it isomorphism
type} of the torsion part of the extended homology
(viewed as an object of the extended abelian category) goes as follows.
First, one refers to theorem of Gromov and Shubin \cite{GS} to obtain that
the dilatational equivalence class of the spectral density function is an 
isomorphism type invariant of a torsion object 
(cf. also \cite{F2}, Proposition 4.5). Secondly, one may
use Lemma 1.20 of \cite{BFKM}, for example, to arrive at the desired 
conclusion. 

Since the extended $L^2$ homology (up to isomorphism) depends 
only on the homotopy type of the chain complex, we obtain 
(cf. also \cite{BFKM}, Proposition 5.6):

\proclaim{3.9. Corollary} The property of a chain complex $C$ of
finitely generated Hilbertian modules to be of the determinant class
depends only on the homotopy type of $C$ in the category of finitely 
generated complexes of Hilbertian modules over $\A$.
\endproclaim

Now we will formulate the generalization mentioned above.

\proclaim{3.10. Proposition} In the category of finite chain complexes
$$C:\quad
(0\to C_N@>{\partial}>> C_{N-1}@>{\partial}>> \dots C_0\to 0)\tag43$$
of finitely generated Hilbertian modules of determinant class 
there is a natural isomorphism
between the determinant lines of graded Hilbertian modules
$$\phi_{C}:\quad \det(C)\to\ \det(H_\ast),\tag44$$
where $H_\ast=\oplus H_i$ is the graded Hilbertian module consisting of
the reduced $L^2$-homology of $C$, that is,
$$H_i=\ker(\partial)/\cl(\im(\partial)).\tag45$$
\endproclaim
\demo{Proof} Consider the following
two sequences:
$$ 0\to \bar B_i\to Z_i\to H_i\to 0$$
and
$$ 0\to Z_i\to C_i\to \bar B_{i-1}\to 0.$$
The first sequence is exact and the second sequence is $D$-exact. By
propositions 2.6 and 3.5 above we have natural isomorphisms
$$\det(Z_i)\to \det(H_i)\otimes\det(\bar B_i)$$
and
$$\det (C_i)\to \det(Z_i)\otimes \det(\bar B_{i-1}).$$
Together they give the natural isomorphism
$$\det (C_i)\to  \det(H_i)\otimes\det(\bar B_i)\otimes \det(\bar B_{i-1})$$
and thus we obtain the natural isomorphism
$$\phi_{C}:\quad \det(C)\ = \ \otimes_i \det(C_i)^{(-1)^i} \ \to \
\ \otimes_i \det(H_i)^{(-1)^i} \  = \det(H_\ast).$$
This completes the proof. \qed
\enddemo

Now we will compute numerically the canonical isomorphism (44).

\proclaim{3.11. Proposition} Suppose that we have fixed an admissible
scalar product on each chain space $C_i$ of the chain complex (43) (which is
assumed to be of determinant class). Let
$\alpha_i\in\det(C_i)$ represent the induced volume form on $C_i$
and let
$$\alpha=\prod\alpha_i^{(-1)^i}\in\det(C_\ast)$$ 
be their alternating product. Let 
$\Delta_i=\partial^\ast\partial+ \partial\partial^\ast: C_i\to C_i$
denote the Laplacian constructed by means of the chosen scalar products.
The reduced homology $H_\ast=H_\ast(C)$ can now be identified with the space
of harmonic forms via the Hodge decomposition.
 We will denote by $\beta\in\det(H_\ast)$ the volume
form on the graded homology inherited by means of the above 
identification. Then the following formula holds
$$\phi_{C} = \prod_{i=0}^N \Det_\tau(\Delta_i^+)^{(-1)^ii/2}\cdot 
\alpha^\ast\otimes \beta  \ \in \det(C)^{-1}\otimes\det(H_\ast),\tag46$$
where $\Delta_i^+$ denotes the restriction of the Laplacian $\Delta_i$
on the orthogonal complement to the space of harmonic forms.
\endproclaim
\demo{Proof} First we observe that the canonical isomorphism (44) behaves
in a very simple way with respect to direct sums of chain complexes. Namely,
$\phi_{C_1\oplus C_2} = \phi_{C_1}\otimes\phi_{C_2}$. 
Now, any chain complex of Hilbertian modules can be represented as a direct
sum of complexes of the form
$$(\dots\to 0\to C_i@>{\partial}>>C_{i-1}\to 0\to \dots )$$
(only two chain spaces are non-zero and the boundary map is a $D$-isomorphism)
and of chain complexes with zero differentials. Using the above remark
one has to check (47) only for such complexes, which is straightforward. \qed
\enddemo
Note, that the analogous formula for {\it cochain complexes} $C$
(of determinant class) is slightly different
$$\phi_{C} = \prod_{i=0}^N \Det_\tau(\Delta_i^+)^{(-1)^{i+1}i/2}\cdot 
\alpha^\ast\otimes \beta  \ \in \det(C)^{-1}\otimes\det(H^\ast),\tag47$$
cf. \cite{BZ}, \cite{BFKM}.

\heading{\bf \S 4. Combinatorial $L^2$-torsion}\endheading

In this section we will define and study a generalization of the classical
construction of the combinatorial torsion (of Reidemeister, Franz and DeRham)
to the case of infinite dimensional representations, which are modules over 
a finite von Neumann algebra. Given a
finite polyhedron $K$ and an unimodular representation of its fundamental
group in a module $M$ over a finite von Neumann algebra $\A$, the torsion
invariant defined here is a {\it positive element of the determinant
line}
$$\det(M)^{-\chi(K)} \otimes \det(H_\ast(K,M)).$$
Under the assumption that the Euler characteristic of $K$ vanishes
(which is always the case if $K$ is a closed manifold of odd dimension),
the torsion does not depend on the choice of the volume form on $M$ and it
can be understood as a volume form on the reduced 
$L^2$-homology $H_\ast(K,M)$, that
is, as an element of the determinant line $\det(H_\ast(K,M))$.

\subheading{4.1} Let $K$ be a finite cell complex.
Denote by $\pi=\pi_1(K)$ its fundamental
group and by $C_\ast(\widetilde K)$ the cellular chain complex of the universal
covering $\widetilde K$ of $K$. Note that the group $\pi$ acts on $C_\ast(\widetilde K)$ 
from the
left and $C_\ast(\widetilde K)$ is a finite complex of free $\Z[\pi]$-modules
with lifts of the cells of $K$ representing a free basis of $C_\ast(\widetilde K)$
over $\Z[\pi]$.

Let $\A$ be a fixed finite von Neumann algebra with a finite, normal and
faithful trace $\tau$. Let $M$ be a finitely generated Hilbertian
module over $\A$. As above, we will denote by $\B(M)=\B_{\A}(M)$ 
the commutant of $M$.

\subheading{4.2} We will consider {\it representations of the group}
$\pi$ in $M$. Any such representation is a homomorphism $\pi\to \B(M)\ =\ \B_{\A}(M)$ 
of multiplicative groups. We will think of $\pi$ as acting
on $M$ from the right (via the representation $\pi\to \B(M)$); thus $M$
will have a structure of $(\A-\pi)$-bimodule. In this situation we will 
say that $M$ is {\it Hilbertian $(\A-\pi)$-bimodule}.

We will say that a Hilbertian $(\A-\pi)$-bimodule $M$ is {\it unimodular} if 
for every element $g\in \pi$
the Fuglede-Kadison determinant
$$\Det_\tau(g)\ =\ 1,\tag48$$
equals 1, where $g$ is viewed as an invertible linear operator 
$M\to M$, given by the right multiplication by $g$.

We will say that a Hilbertian $(\A-\pi)$-bimodule $M$ is {\it unitary} if 
there is exists an admissible
scalar product $\<\ ,\ \>$ on $M$ such that the action of $\pi$ preserves this
scalar product. Obviously, any unitary Hilbertian bimodule is unimodular.

Any Hilbertian $(\A-\pi)$-bimodule $M$ determines a {\it flat Hilbertian
bundle over $K$ with fiber $M$}. In fact, consider Borel's construction 
$$\M\ =\ M\times_{\pi}{\widetilde K}$$ 
together with the obvious projection map $\M\to K$; it has a canonical
structure of a flat $\A$-Hilbertian bundle. Equivalently, it can be viewed as
a locally
free sheaf of Hilbertian $\A$-modules.

\subheading{4.3. Examples} As a first concrete example, consider 
the classical case, when 
$\A={\Cal N}(\pi)$ is the von Neumann algebra of $\pi$ and
$M=\ell^2(\pi)$ is the completion of the group algebra of $\pi$
with respect to the canonical trace on it, with $\A$ acting on $M$ 
from the left and with $\pi$ acting
on $M$ from the right. This $(\A-\pi)$-bimodule $M=\ell^2(\pi)$ is
unitary and thus unimodular, as well.

As a more general example consider the following Hilbertian 
$(\A-\pi)$-bimodule $M\ =\ \ell^2(\pi)\otimes_{\C}V$ where $V$ is a finite
dimensional unimodular representation of $\pi$. Here the left action of 
$\A={\Cal N}(\pi)$ on $M$ is the same as the action on $\ell^2(\pi)$ and
the right action of $\pi$ is the diagonal action:
$(x\otimes v)g=xg\otimes vg$ for $x\in \ell^2(\pi)$, $v\in V$, and $g\in \pi$.

\subheading{4.4. Construction of the torsion} Given an unimodular
$(\A-\pi)$-bimodule $M$, we can
form the complex
$$C_\ast(K,M)\ =\ M\otimes_{\Z\pi}C_\ast(\widetilde K).\tag49$$
It is a chain complex of finitely generated $\A$-modules; each chain group
$C_i(K,M)$ can be identified with a finite direct sum of a number of copies
of $M$ and the number of summands equals the number of $i$-dimensional
simplexes of $K$.

We will assume that {\it this chain complex $C_\ast(K,M)$ is 
of determinant class}, cf. 3.8 above.
Then by Proposition 3.10 we obtain a natural isomorphism
$$\det(C_\ast(K,M))@>{\simeq}>>\det(H_\ast(K,M)),\tag50$$
where $H_\ast(K,M)$ denotes {\it the reduced $L^2$-homology} of $C_\ast(K,M)$. 
Now, under the unimodularity
assumption there is natural isomorphism
$$\det(M)^{\chi(K)}@>{\simeq}>> \det(C_\ast(K,M)),\tag51$$
defined as follows.
For any cell $e\subset K$ fix a lifting $\tilde e$ of $e$ in the universal
covering. Then the cells $\tilde e$ form a free $\Z\pi$-basis of the
complex $C_\ast(\widetilde K)$ and therefore they allow us to represent $C_\ast(K,M)$
as the direct sum of copies of $M$, one for each cell. Thus, the determinant
line $\det(C_\ast(K,M))$ can be identified with $\det(M)^{\chi(K)}$,
since the cells of odd dimension contribute negative factors of $\det(M)$
into the total determinant line.

We only have to show that this identification
does not depend on the choice of the liftings $\tilde e$. Consider an arbitrary
set of liftings of the cells of $K$. It can be obtained as follows.
Suppose that for each cell $e\subset K$ an element $g_e\in\pi$ has been fixed.
Then the cells $g_e\tilde e$ form another set of liftings. 
Since the Fuglede-Kadison determinant of the map
$$\oplus M\ \to \oplus M\tag52$$
(where in the sums on both sides the number of copies of $M$ is equal
to the number of the cells in $K$), given by the diagonal matrix with
$g_e$ on the diagonal, is 1 (since the representation is unimodular),
we see that the isomorphism (51)
is canonical.

\subheading{Definition} One may interpret the composition of isomorphisms (51) 
and (50) as a nonzero element of the line
$$\rho_K\ =\ \rho_{K,M}\ \ \in\ \det(M)^{-\chi(K)}\otimes \det(H_\ast(K,M)),\tag53$$   
which will be called {\it combinatorial $L^2$ torsion, or $L^2$ 
Reidemeister - Franz torsion}.

In \S 5 we will consider analytic version of this construction.

\subheading{4.5. Remarks} 

1. Note that in the case $\A=\C$ we arrive to the
classical definitions, cf. \cite{Mi1}, \cite{Mi2}, \cite{Mi3}, 
\cite{Mu}, \cite{BZ}.

2. Recall that the classical Reidemeister-Franz torsion is not, in general,
a homotopy invariant,
so one cannot expect homotopy invariance from our torsion invariant.
But the {\it combinatorial invariance} holds, cf. below.

3. Although our notation $\rho_K$ for the torsion invariant does not
involve explicitly the trace $\tau:\A\to \C$, the whole construction
(including the Fuglede-Kadison determinants and the determinant lines)
certainly depend of the choice of the trace $\tau$. Thus, in fact,
we have one determinant line
$\det(M)^{-\chi(K)}\otimes\det(H_\ast(K,M))$
for each trace $\tau$, forming a bundle, and $\rho_K$ is a section of this
bundle. 

Consider, for example, the case when the initial von Neumann algebra $\A$
is the group ring of a finite group, $\A=\C[G]$. Then we have one trace for
any irreducible representation of $G$. Thus, we have one determinant line
for each irreducible representation, and the torsion $\rho_K$ is a function
on the classes of the irreducible representations with values in these
determinant lines.

4. In the case when $\chi(K)=0$, the torsion invariant is just an element
of the determinant line of the reduced $L^2$-homology
$\det(H_\ast(K,M))$. This is always the case if
$K$ is a closed manifold of odd dimension.

5. If the representation $M$ is {\it unitary}, then the scalar product on
$M$ determines a well defined element in $\det(M)$. Then, again, we can
consider the torsion as a volume on the reduced $L^2$-homology, that is,
as an element of $\det(H_\ast(K,M))$.

Assuming additionally, that the reduced $L^2$-homology $H_\ast(K,M)$ vanishes, 
we can identify the determinant line $\det(H_\ast(K,M))$ with $\R$ and
so $\rho_K$ is just a number. Under this assumption it was studied in
\cite{CM}, \cite{Lu}, \cite{LR}.

\proclaim{4.6. Theorem (Combinatorial Invariance)} Let $K$ be a finite
polyhedral cell complex and let
$K^\prime$ be its subdivision. Suppose that $M$ is an unimodular
representation of $\pi=\pi_1(K)$ over a finite von Neumann algebra $\A$
with a trace $\tau$ and suppose the the complex $M\otimes_{\Z\pi}C_\ast(\widetilde K)$
is of determinant class. Let
$$\psi:H_\ast(K,M)\to H_\ast(K^\prime,M)\tag54$$
be the isomorphism induced on the reduced $L^2$ homology by the 
subdivision chain map.
Then $\psi$ is an isomorphism of Hilbertian modules and
the induced by $\psi$ map
$$\id\otimes\psi_\ast:
\ \det(M)^{-\chi(K)}\otimes\det(H_\ast(K,M))\ \to\
\det(M)^{-\chi(K)}\otimes\det(H_\ast(K',M))\tag55$$
maps $\rho_{K}$ onto $\rho_{K'}$.
\endproclaim
\demo{Proof} It is enough to consider the elementary subdivision when a
single $q$-dimensio\-nal cell $e$ is divided into two $q$-dimensional cells
$e_+$ and $e_-$ introducing an additional separating $(q-1)$-dimensional
cell $e_0$, see the figure.

\midspace{4cm}\caption{Figure 1}

We have the exact sequence of free left $\Z[\pi]$-chain complexes
$$0\to C_\ast(\widetilde K)@>{\psi}>>C_\ast(\widetilde K^\prime)\to 
D_\ast\to 0\tag56$$
where the chain complex $D_\ast$ has nontrivial chains only in dimension
$q$ and $q-1$ and $D_q$ and $D_{q-1}$ are both free of rank one. The free
generator of the module $D_q$ can be labeled with $e_+$ and the generator
of $D_{q-1}$ can be labeled with the cell $e_0$ and then the boundary 
homomorphism is given by $\partial(e_+)=e_0$.

Fix an admissible scalar product $\<\ ,\ \>$ on $M$. 
Using the cell structures of $K$ and $K^\prime$, we then
obtain canonically the scalar products on the complexes
$$M\otimes_{\Z\pi} C_\ast(\widetilde K),\quad    
M\otimes_{\Z\pi} C_\ast(\widetilde K^\prime),\quad  
M\otimes_{\Z\pi} D_\ast,$$  
and, according to our definitions, these scalar products represent volume
forms
$$x\in \det(M\otimes_{\Z\pi} C_\ast(\widetilde K)),\quad 
y\in\det(M\otimes_{\Z\pi} C_\ast(\widetilde K^\prime)),\quad   
z\in\det(M\otimes_{\Z\pi} D_\ast)$$  

Note that the chain complex $D_\ast$ is acyclic and so the canonical 
isomorphism
(44), which we will denote in this case by $\phi_{D_\ast}$, identifies the 
determinant line of $D_\ast$ with $\R$. It is easy to see that
$\phi_{D_\ast}(z)\ =\ 1.$

Using the exact sequence (56), we obtain the commutative diagram
$$
\CD
\det(M\otimes_{\Z\pi} C_\ast(\widetilde K^\prime))@>>>\det(M\otimes_{\Z\pi} C_\ast(\widetilde K))\otimes
\det(M\otimes_{\Z\pi} D_\ast)\\
@V{\phi_{K^\prime}}VV    @VV{\phi_K\otimes \phi_{D_\ast}}V\\
\det(H_\ast(K^\prime, M))@>>{\psi_\ast^{-1}}>\det(H_\ast(K,M))
\endCD
$$
and the upper horizontal map maps $y$ into $x\otimes z$. By the definition,
we have 
$$\rho_{K^\prime}\ =\ \<\ ,\ \>^{-\chi(K)}\otimes \phi_{K^\prime}(y), \quad
\text{and}\quad
\rho_{K}\ =\ \<\ ,\ \>^{-\chi(K)}\otimes \phi_{K}(x) $$
and so from the diagram above we get 
$\id\otimes \psi_\ast^{-1}(\rho_{K^\prime})=\rho_K$. $\square$
\enddemo
\subheading{4.7. Cohomological version of the construction} 
Here we will mention a variation of the construction of the combinatorial
torsion (53), based on cohomology instead of homology. 

Suppose that
$M$ is {\it a Hilbertian $(\pi - \A)$-bimodule}. This notion is similar to
the notion of a Hilbertian $(\A - \pi)$-bimodule, introduced in 4.2. It means
that $M$ is a topological vector space supplied with a right $\A$-action
and a left $\pi$-action which commute with each other and such that $M$
is a finitely generated Hilbertian module over $\A$.

If $K$ is a finite polyhedron, consider the cellular chain complex 
$C_\ast(\widetilde K)$ of its universal covering $\widetilde K$. Then 
$$C^\ast(K,M)\ =\ \overline{\Hom_{\Z\pi}(C_\ast(\widetilde K),M)}$$
is a cochain complex of finitely generated left Hilbertian modules
over $\A$. Here the bar means that we convert the right $\A$-module 
structure on $\Hom_{\Z\pi}(C_\ast(\widetilde K),M)$ into a left structure
using the involution of $\A$. The reduced cohomology of $C^\ast(K,M)$
will be denoted 
$$H^\ast(K,M).\tag57$$

We suppose now that {\it the chain complex $C^\ast(K,M)$ is of determinant class
and that the $(\pi-\A)$-module $M$ is unimodular} (this means that for
each $g\in \pi$ the map $M\to M$, given by multiplication on $g$ from the 
left has Fuglede-Kadison determinant 1). Then there is a natural isomorphism 
$$\det(M)^{\chi(K)}\ \to \ \det(C^\ast(K,M)),$$
which is similar to (51). Composed with the isomorphism of Proposition 3.10
$$\det(C^\ast(K,M))\ \to\ \det(H^\ast(K,M))$$
(where $H^\ast(K,M)$ denotes the reduced $L^2$ cohomology) 
it gives a canonical element
$$\rho_K\ =\ \rho_{K,M}\ \in \det(M)^{-\chi(K)}\otimes\det(H^\ast(K,M)).
\tag58$$
It will be called {\it combinatorial $L^2$ torsion, or $L^2$ Reidemeister - 
Franz torsion}.

Again, similar to 4.5, remark 5, we may view the torsion $\rho_K$ as an element
of the line $\det(H^\ast(K,M))$, assuming that the representation $\pi$ is unitary
and so $M$ has a specified admissible invariant scalar product.

\heading{\bf \S 5. Analytic $L^2$-torsion}\endheading

In this section we
generalize the classical construction of analytic torsion
of D.B. Ray and I.M. Singer \cite{RS}
to the case of infinite dimensional representations of the fundamental group.  
The invariant we construct represents
a volume form on the reduced $L^2$ cohomology. We will see in the next section 
that this allows really understand the meaning of
the recent theorem of Burghelea, Friedlander,
Kappeler and McDonald \cite{BFKM}.

Throughout this section $\A$ denotes a finite von Neumann algebra 
with a fixed finite, normal, and faithful trace $\tau$. 

\subheading{5.1 Flat Hilbertian $\A$-Bundles} Let $M$ be a finitely generated 
Hilbertian $(\pi - \A)$-bimodule (cf. 4.2). Recall that this means that 
$\A$ acts on $M$ 
from the right and $M$, with respect to this action, is a finitely generated 
Hilbertian module (cf. 1.3); also, $\pi$ is a discrete group acting on $M$
from the left and the action of $\pi$ commutes with that of $\A$. 

Let $X$ be a connected, closed, smooth 
manifold with fundamental group $\pi$. 
{\it A flat Hilbertian $\A$-bundle with fiber $M$ 
over $X$} is an associated bundle $p:{\E}\to X$, where 
$\E\ =\ (M\times\widetilde X)/\sim$ with its natural projection onto $X$.
Here $(v,x) \sim (gv,gx)$ for all $g\in \pi$, $x\in \widetilde X$ and 
$v\in M$.
Here $\widetilde X$ denotes the universal covering of $X$. The map
$p:{\E}\to X$ is then
a locally trivial bundle of topological vector spaces (cf. chapter 3 of 
\cite{La}),
which has a natural fiberwise left action of $\A$.

Any smooth section $s$ of $\E\to X$ can be uniquely represented by
a smooth equivariant map $\phi:\widetilde X\to M$, that is, $\phi(gx) 
=g\phi(x)$ for all $g\in\pi$ and $x\in \widetilde X$. 

\subheading{5.2} Given a flat Hilbertian 
$\A$-bundle $\E\to X$ over a closed connected manifold $X$,
one can consider the space of smooth differential $j$-forms on $X$ 
with values in $\E$. This space will be denoted by $\Omega^j(X,\E)$. 
It is naturally 
defined as a left $\A$-module. An element of $\Omega^j(X,\E)$
can be uniquely represented as a $\pi$-invariant differential form in
$M\otimes_\C \Omega^j(\widetilde X)$. Here one considers the total (diagonal) $\pi$ action,
that is, the tensor product of the actions of $\pi$ on $\Omega^j(\widetilde X)$
and on $M$. More precisely, if $\omega\in\Omega^j(\widetilde X)$ and $v\in M$, then
$v\otimes\omega$ is said to be $\pi$-invariant if $vg^{-1}\otimes g^*\omega 
= v\otimes\omega$ for all $g\in\pi$.

{\it A flat $\A$-linear connection}
on a flat Hilbertian $\A$-bundle $\E$ is defined as an $\A$-homomorphism
$$\nabla : \Omega^j(X,\E)\to \Omega^{j+1}(X,\E)$$
such that 
$$\nabla(f\omega ) = df\wedge\omega + f\nabla (\omega) 
\qquad\text{and}\qquad\nabla^2 = 0
$$
for any $\A$ valued function $f$ on $X$ and 
for any $\omega\in \Omega^j(X,\E)$.
On the flat Hilbertian $\A$-bundle $\E$, as defined in the previous paragraph,
there is a {\it canonical flat
$\A$-linear connection} $\nabla$, which is given as follows: under the
identification of $\Omega^j(X,\E)$ given in the previous paragraph, 
one defines $\nabla(v\otimes \omega) = v\otimes d\omega$, where $d$ is the 
De Rham exterior derivative. 

\subheading{5.3 Hermitian metrics and L$^2$ scalar products}
{\it A Hermitian metric $h$} on a flat Hilbertian $\A$-bundle $p:\E\to X$
is a smooth family of admissible scalar products on the fibers. Any
Hermitian metric on $p:\E\to X$ defines a wedge product 
$$
\wedge : \Omega^i(X,\E)\otimes \Omega^{j}(X,\E)\rightarrow
\Omega^{i+j}(X)
$$
similar to the finite dimensional case. 

A Hermitian metric on $p:\E\to X$ together with a Riemannian metric
on $X$ determine an scalar product on $\Omega^i(X,\E)$ in the
standard way; namely, using the Hodge star operator 
$$\ast:\Omega^j(X,\E)\to \Omega^{n-j}(X,\E)$$
one sets
$$(\omega,\omega^\prime)\ =\ \int_X \omega\wedge\ast\omega^\prime. $$

With this scalar product $\Omega^i(X,\E)$
becomes a pre-Hilbert space. Define the space of
$L^2$ differential $j$-forms on $X$ with coefficients in
$\E$, which is denoted by $\Omega_{(2)}^j(X,\E)$, to be the Hilbert
space completion of $\Omega^j(X,\E)$. 

\subheading{5.4 Reduced L$^2$ cohomology} 
Given a flat Hilbertian $\A$ bundle $p:\E\to X$, one 
defines the {\it reduced $L^2$ cohomology of $\E$} as the quotient 
$$
H^{j}(X,\E)=\frac{\ker \nabla/\Omega^{j}_{(2)}(X,\E)}{\cl(\im\;\nabla/
\Omega^{j-1}_{(2)}(X,\E))},\tag59
$$
where the connection $\nabla$ on $\E$ is assumed to be extended
to an unbounded, densely defined operator $\Omega_{(2)}^j(X,\M)\to 
\Omega_{(2)}^{j+1}(X,\M)$.
Then $H^j(X,\E)$ is naturally defined as a Hilbertian module over $\A$.
Considered as Hilbertian 
module over $\A$, it is independent on the choice of the Riemannian metric
on $X$ and the Hermitian metric on $\E$.

The corresponding $L^2$ Betti numbers are denoted by 
$$
b^{j}(X,\M) = \dim_\tau \left( H^{j}(X,\M)\right).\tag60
$$

Given a smooth triangulation of the manifold $X$, one may consider also the
reduced $L^2$-cohomology defined combinatorially (57) with coefficients in 
the monodromy representation $M$ of $\E$. The De Rham type
theorem states that these two kinds of reduced $L^2$ cohomology
$H^\ast(X,M)$ (the combinatorial (57)) and $H^\ast(X,\E)$ (the analytic (59))
{\it are canonically isomorphic viewed as Hilbertian $\A$-modules}. 
The canonical isomorphism 
$$I :  H^\ast(X,\E) \to  H^\ast(X,M)\tag61$$
is given by integration of forms along the simplices. The isomorphism (61) was
established by J. Dodziuk in \cite{D}. In fact, Dodziuk studied only the case of the
regular representation but his arguments can be easily generalized to the
case of an arbitrary flat bundle.

\subheading{5.5} Let 
$\Delta_j = \int_0^\infty \lambda dE_j(\lambda)$ denote the spectral
decomposition of the Laplacian
$\Delta_j=\nabla \nabla^{*} + \nabla^* \nabla:\Omega_{(2)}^j(X,\M)
\rightarrow\Omega_{(2)}^j(X,\M),$
where $\nabla^{*}$ denotes the formal adjoint of $\nabla$ with respect to the
$L^{2}$ scalar product on $\Omega_{(2)}^j(X,\M)$. 
Note that by definition, the
Laplacian is a formally self-adjoint operator, which is densely defined. We 
also denote by $\Delta_j$ the self adjoint extension of the Laplacian.

 Recall, that {\it the spectral density 
function} is defined as
$N_j(\lambda) = \Tr_\tau(E_j(\lambda))$ and the {\it theta function} 
is by definition
$$\theta_j(t) = \int_0^\infty e^{-t\lambda} dN_j(\lambda)= 
\Tr_\tau(e^{-t\Delta_j}) -
b^{j}(X,\M).$$ 
Here we use the well known fact that the projection
$E_j(\lambda)$ and the heat operator $e^{-t\Delta_j}$ have smooth 
Schwartz kernels, which are smooth sections of a bundle over $X\times X$ with
fiber the commutant of $M$, cf. \cite{BFKM}, \cite{GS}, \cite{Luk}. 
The symbol $\Tr_\tau$ denotes 
application of the canonical trace
on the commutant (cf. Proposition 1.8) to the restriction of the kernels to 
the diagonal followed by integration over the manifold $X$. See also \cite{M},
\cite{L}  and \cite{GS} for the case of the flat bundle defined by the 
regular representation of the fundamental group.

\subheading{5.6 Definition} A flat Hilbertian $\A$-bundle $\M\rightarrow X$ 
is said to be of {\it determinant 
class} if
$$
   \int_0^1 \log (\lambda) dN_j(\lambda) >\ -\infty\qquad
   \text{or, equivalently}\quad \int_1^\infty t^{-1} \theta_j(t) dt < \infty.  \tag62
$$
for all $j=0,....,n$. Equivalence of these two conditions was proved in
\cite{BFKM}, Proposition 2.12.

This notion has been studied earlier. 
In \cite{CM}, this property was called $D$-acycl\-icity.
In \cite{BFKM}, the term {\it a-determinant class} was used.
In \cite{BFKM}, D. Burghe\-lea et al. introduced also the notion of 
{\it flat bundles of c-determinant class} (here "a" stands for {\it analytic} 
and "c"
stands for {\it combinatorial}); these two notions are actually equivalent
as shown by A.V.Efremov \cite{E}. The c-determinant class condition
is
obviously equivalent to Definition 3.8. This shows in particular that the 
{\it determinant class property of a flat Hilbertian $\A$ bundle does
not depend on the choice of metrics $g$ on $X$ and $h$ on the bundle $\E$}. 

\subheading{5.7. Definition (\cite{BFKM})}
For $\lambda >0$, the {\it zeta function of the Laplacian} $\Delta_j$ is
defined
on the half-plane $\Re(s)>n/2$ as
$$
\zeta_{j}(s, \lambda, \E)  = \frac{1}{\Gamma(s)}
\int^{\infty}_{0}t^{s-1} e^{-\lambda t}\theta_{j}(t)dt.\tag63
$$

It has been shown by \cite{BKFM}, that $\zeta_{j}(s, \lambda, \E)$, 
viewed as a function of $s$, is holomorphic
in the half-plane $\Re(s)>n/2$ (where $n=\dim X$)
and has a meromorphic continuation to $\C$ with no pole at $s=0$.
Under the hypothesis that the flat
Hilbertian $\A$-bundle $\E\to X$ is of determinant class,
the limit $\lim_{\lambda\to 0}\zeta_{j}'(0, \lambda, \E)$ exists, where
the prime denotes the differentiation with respect to $s$; we will denote
this limit by $\zeta_j'(0,0,\E)$. Denote also 
$$\zeta'(0,0,\E) = \sum_{j}(-1)^jj\zeta_j'(0,0,\E).\tag64$$

\subheading{5.8. The construction of $L^2$ analytic torsion } 
Let $(X,g)$ be a closed 
Riemannian manifold of dimension $n$ with $\pi=\pi_{1}(X)$ and 
let $\E$ be a flat Hilbertian $\A$-bundle over $X$ with a Hermitian metric $h$ 
and with fibre $M$. We assume that $\E$ is of determinant class.

Using the Hodge theorem (cf. \cite{D}, \cite{BFKM}, \cite{GS}) one obtains 
an identification between the reduced $L^2$ cohomology $H^{j}(X,\E)$ (defined
by (59)) and the space of harmonic forms , i.e. the kernel of 
the Laplacian $\Delta_j$ acting on $\Omega_{(2)}^j(X,\M)$. Since the last space
is embedded into a Hilbert space it inherits an admissible scalar product and thus
we obtain an admissible scalar product on the cohomology space $H^{j}(X,\E)$.

These admissible scalar products on $H^{j}(X,\E)$, 
determine nonzero elements of the determinant lines $\det(H^{j}(X,\M))$
for all $j$ and thus, their alternating product in the line
$$\det(H^{\ast}(X,\M))
=\prod_{j=0}^n \det(H^{j}(X,\M))^{(-1)^j}\tag65$$
is defined; the last element we will denote by $\tilde\rho(g, h)$. 
This notation emphasizes dependence on the metrics $g$ and $h$.

\subheading{ Definition} {\it The $L^{2}$ analytic torsion} 
is defined as the element of the determinant line
$\rho_{\E}(g,h)\in\det(H^{\ast}(X,\M))$, where
$$\rho_{\E}(g, h)=e^{\frac{1}{2} \zeta'(0,0,\M)}\cdot\tilde\rho(g, h).\tag66$$
Here $\zeta'(0,0,\M)$ is defined in (64). Observe that the sign in the exponent
of (66) is different from the one in the 
usual formulas (as in \cite{BZ}, \cite{BGS}); it happens because we consider the
elements of the corresponding lines (the volume forms)
but not the metrics on them.

Thus, the $L^2$ analytic torsion represents a volume form on the reduced 
$L^2$ cohomology. We will see later in section 6, that $\rho_{\E}(g, h)$
does not depend on the metrics $g$ and $h$ if $\dim(X)$ is odd.

Note, that in the case when ${\Cal A}=\C$, we
arrive at the classical definition of the Ray-Singer-Quillen
metric on the determinant of the cohomology.

Assuming that the reduced $L^{2}$ cohomology $H^*(X,\M)$
vanishes, we can identify canonically the determinant line 
$\det({H}^{*} (X,\M))$ with $\R$, and so the torsion $\rho_{\E}$ just defined
turns into a number. In this case, it  coincides with the invariant studied
in \cite{M} and \cite{L}. 

\heading{\bf \S 6. Theorem of Burghelea,  Friedlander, Kappeler and McDonald}
\endheading

In this section we will show that a recent theorem of Burghelea,  Friedlander, 
Kappeler and McDonald \cite{BFKM} can be naturally formulated in terms of the 
notions introduced in the present paper. Namely, we will show that the main 
theorem of \cite{BFKM} (Theorem 2) 
is equivalent to the statement that the $L^2$ combinatorial and $L^2$
analytic torsion invariants, as defined in the present paper,
determine identical volume forms on the reduced $L^2$
cohomology (cf. Theorem 6.1 below for the precise statement). Note that the
original theorem of Burghelea et al. \cite{BFKM} 
is formulated using three numerical torsion type invariants,
neither of which is a topological invariant.

We are very thankful to the referee, who pointed out this statement to us.

\proclaim{6.1. Theorem} Let $X$ be a closed odd-dimensional Riemannian
manifold
and let $\E\to X$ be a flat Hilbertian bundle (cf. 5.1) of determinant class
supplied with a flat fiberwise Hermitian metric $h$. Consider the reduced $L^2$
cohomology $H^\ast(X,\E)$ and the analytic $L^2$ torsion 
$\rho_{\E}(g,h)\in\det(H^{\ast}(X,\M))$. On the other hand, consider a smooth
triangulation of $X$ and the combinatorially defined reduced $L^2$ cohomology
$H^\ast(X,M)$ (cf. (57), where $M$ denotes the fiber of $\E$, viewed together 
with the monodromy representation of the fundamental group of $X$ and with the
induced admissible scalar product. Then the $L^2$ combinatorial torsion
$\rho_X\in \det(H^\ast(X,M)$ is defined, cf. (58).  Then the isomorphism 
between the cohomological determinant lines
$I_\ast:  \det(H^\ast(X,\E)) \to \det(H^\ast(X,M))$ determined by the
De Rham isomorphism (61) (cf. 2.3) maps the analytically defined torsion
$\rho_\E(g,h)$ onto the combinatorially defined torsion $\rho_X$, i.e.
$$I_\ast(\rho_{\E}(g,h)) = \rho_X.\tag67$$
\endproclaim

In particular we obtain that the analytic torsion $\rho_{\E}(g,h)$ does not depend
on the Riemannian metric $g$ on $X$ and Hermitian metric $h$ on $\E$.

\demo{Proof} The proof consists of interpreting the result of \cite{BFKM}
in terms of the notions of this paper. 
Let us denote by $\alpha\in\det(H^\ast(X,M))$ the volume form on the 
combinatorially defined reduced $L^2$ cohomology correspondidng to the given
triangulation (i.e. corresponding to the metric on the combinatorial
harmonic forms via the Hodge isomorphism). Denote by 
$\beta\in\det(H^\ast(X,\E))$ the volume form on the analytically defined 
reduced $L^2$ cohomology which corresponds to the induced admissible
scalar product on the space of harmonic forms via the Hodge isomorphism.

We claim that
$$I_\ast(\beta) = T_{met}^{-1}\cdot\alpha,\tag68$$
where $T_{met}$ ({\it the torsion of the metric}) is a numerical invariant,
measuring the relative size of the Riemannian metric and the triangulation, and
is defined in the introduction to \cite{BFKM}. Formula (68) follows by 
comparing the definition of \cite{BFKM} with our definitions: formulae 
(26) and (23).

Now, again by comparing the definitions of \cite{BFKM} with ours,
we obtain for the combinatorial torsion
$$\rho_X = T_{comb}\cdot\alpha,\tag69$$
and also
$$\rho_\E(g,h) = T_{an}\cdot \beta\tag70$$
for the analytic torsion, where $T_{comb}$ and $T_{an}$ are the 
numbers introduced in \cite{BFKM}.

Therefore, combining (68), (69) and (70), we obtain
$$
\aligned
I_\ast(\rho_\E(g,h)) &= T_{an}\cdot I_\ast(\beta) =\\
&= T_{an}\cdot T_{met}^{-1}\cdot\alpha =\\
&= T_{comb}\cdot\alpha = \rho_X. 
\endaligned
$$
This completes the proof. \qed
\enddemo

Note that Theorem 6.1 implies Theorem 4.6 in the case when the polyhedron is
a manifold of odd dimension and the bundle is unitary.

It is plausible that Theorem 6.1 holds without assuming that the flat 
bundle $\E\to X$ is unitary.  In fact, the analytic torsion (66) is defined without
assumptions of this kind. One may show that $\rho(g,h)$
is independent of the metrics
$g$ and $h$, if the dimension of $X$ is odd (this is a general property of 
analytic torsion for elliptic complexes, cf. \cite{F3}). In order to be able to
remove the assumption of unitarity (or unimodularity)  from Theorem 6.1 
one has to be able first to
define the combinatorial torsion (58) without assuming unimodularity of the bundle
$\E\to X$. 
It is well known that it is impossible for polyhedra;  but it is possible, if one
restricts to orientable manifolds $X$ of odd dimension. A purely combinatorial
construction (using Poincar\'e duality) was suggested in 
\cite{F} for finite dimensional flat bundles; the corresponding metric on the
determinant line of the cohomology was called {\it Poincar\'e - Reidemeister
metric}. We conjecture that the construction of the Poincar\'e - Reidemeister 
metric of \cite{F} 
can be generalized to flat Hilbertian bundles of determinant class and 
that Theorem 6.1 can then be generalized as stating that the 
Poincar\'e - Reidemeister norm (defined combinatorially)
of the analytic torsion (66) equals 1. Compare  \cite{F}, where a similar
theorem was proven (using theorem of Bismut and Zhang \cite{BZ})
in the finite dimensional case.

\enddocument